\documentclass[10pt]{article}
\usepackage{amstext}
\usepackage{amsfonts}
\usepackage{amssymb}
\usepackage{amsbsy}
\usepackage{latexsym}
\usepackage{xy}
\usepackage{hhline}
\xyoption{all}
\newcommand{\vtx}[1]{*+[o][F-]{\scriptscriptstyle #1}} 
\newcounter{num}[section] %

\newenvironment{theo}
{\refstepcounter{num}%
\bigskip\noindent{\bf Theorem~\arabic{section}.\arabic{num}. }\it}

\newenvironment{cor}
{\refstepcounter{num}%
\bigskip\noindent{\bf Corollary~\arabic{section}.\arabic{num}. }\it}

\newenvironment{lemma}
{\refstepcounter{num}%
\bigskip\noindent{\bf Lemma~\arabic{section}.\arabic{num}. }\it}

\newcommand{\example}
{\refstepcounter{num}%
\bigskip\noindent{\bf Example~\arabic{section}.\arabic{num}.}}

\newcommand{\remark}
{\refstepcounter{num}%
\bigskip\noindent{\bf Remark~\arabic{section}.\arabic{num}.}}

\newcommand{\definition}[1]
{\refstepcounter{num}%
\bigskip\noindent{\bf Definition~\arabic{section}.\arabic{num}}~({\it #1}).}

\newenvironment{eq}{\begin{equation}}{\end{equation}}

\newcommand{\Ref}[1]{(\ref{#1})}

\newcommand{\si}{\sigma}
\newcommand{\al}{\alpha}
\newcommand{\be}{\beta}
\newcommand{\ga}{\gamma}
\newcommand{\la}{\lambda}

\newcommand{\ov}[1]{\overline{#1}}
\newcommand{\un}[1]{{\underline{#1}} }

\newcommand{\tr}{\mathop{\rm tr}}

\newcommand{\mdeg}{\mathop{\rm mdeg}}
\newcommand{\diag}{\mathop{\rm diag}}
\newcommand{\Char}{\mathop{\rm char}}

\newcommand{\algA}{\mathcal{A}}    

\newcommand{\M}{\mathcal{M}} 
\newcommand{\N}{\mathcal{N}} 
\newcommand{\I}{{\mathcal I}} 

\newcommand{\FF}{{\mathbb{F}}}   
\newcommand{\NN}{{\mathbb{N}}}
\newcommand{\ZZ}{{\mathbb{Z}}}   

\newcommand{\DP}{{\rm DP} }

\newcommand{\Q}{\mathcal{Q}}    

\newcommand{\Sp}{S\!p}

\newcommand{\loopR}[3]{%
\begin{picture}(20,0)(#1,#2)
\put(-2,1){\llap{$\scriptstyle #3$}} \put(11,3){\circle{20}} \put(20,6){\vector(1,-4){1}}
\end{picture}}
\newcommand{\loopL}[3]{%
\begin{picture}(20,0)(#1,#2)
\put(22,1){$\scriptstyle #3$} \put(9,3){\circle{20}} \put(0,6){\vector(-1,-4){1}}
\end{picture}}

\begin{document}
\title{Orthogonal matrix invariants}
 \author{
A.A. Lopatin \\
Institute of Mathematics, \\
Siberian Branch of \\
the Russian Academy of Sciences, \\
Pevtsova street, 13,\\
Omsk 644099 Russia \\
artem\underline{ }lopatin@yahoo.com \\
}
\date{} 
\maketitle

\begin{abstract} 
The orthogonal group acts on the space of several $n\times n$ matrices by simultaneous conjugation.  For an infinite field of characteristic different from two, relations between generators for the algebra of invariants are described. As an application, the maximal degree of elements of a minimal system of generators is described with deviation $3$. 

This note contains concise but precise description of the results. All proofs can be found in arXiv: 0902.4266 and arXiv: 1011.5201. 
\end{abstract}



\section{Introduction}\label{section_intro}


All vector spaces, algebras, and modules are over an infinite field $\FF$.  By an algebra we always mean an associative algebra.

Let a linear group $G$ be a subgroup of $GL(n)$ and  
$$V=\FF^{n\times n}\oplus\cdots\oplus \FF^{n\times n}$$
be $d$-tuple of $n\times n$ matrices over $\FF$. The group $G$ acts on $V$ by the
diagonal conjugation, i.e.,
\begin{eq}\label{eq_diag_conj}
g\cdot (A_1,\ldots,A_d)=(g A_1 g^{-1},\ldots,g A_d g^{-1}),
\end{eq}%
where $g\in G$ and $A_1,\ldots,A_d\in \FF^{n\times n}$.

The coordinate ring of $V$ (i.e.~the ring of polynomial functions $f:V\to \FF$) is
the ring of polynomials  %
$$\FF[V]=\FF[x_{ij}(k)\,|\,1\leq i,j\leq n,\, 1\leq k\leq d],$$ %
where $x_{ij}(k)$ stands for the coordinate function on $V$ that takes 
$(A_1,\ldots,A_d)\in V$ to the $(i,j)^{\rm th}$ entry of the matrix $A_k$. Denote by
$$X_k=\left(\begin{array}{ccc}
x_{11}(k) & \cdots & x_{1n}(k)\\
\vdots & & \vdots \\
x_{n1}(k) & \cdots & x_{nn}(k)\\
\end{array}
\right)$$%
the $k^{\rm th}\!\!$ {\it generic} matrix ($1\leq k\leq n$).

The action of $G$ on $V$ induces the action on $\FF[V]$ as follows: $(g\cdot
f)(h)=f(g^{-1}\cdot h)$ for all $g\in G$, $f\in \FF[V]$,
$h\in V$. In other words,  %
$$g\cdot x_{ij}(k)= (i,j)^{\rm th}\text{ entry of }g^{-1}X_k g.$$%
The algebra of {\it matrix $G$-invariants} is
$$\FF[V]^{G}=\{f\in \FF[V]\,|\,g\cdot f=f\;{\rm for\; all}\;g\in G\}.$$

If $G$ is a classical linear group, i.e., $G$ belongs to the list $GL(n)$, $O(n)$, $\Sp(n)$, $SO(n)$, $SL(n)$, then a generating set for the algebra $\FF[V]^G$ is known for an arbitrary characteristic of $\FF$ (see~\cite{Donkin92a},~\cite{Zubkov99},~\cite{Lopatin_so_inv}), where we assume that  $\Char{\FF}\neq2$ in the case of $O(n)$ and $SO(n)$. 

Denote coefficients in the characteristic polynomial
of an $n\times n$ matrix $X$ by $\si_t(X)$, i.e., %
\begin{eq}\label{eq1_intro} 
\det(X+\la E)=\sum_{t=0}^{n} \la^{n-t}\si_t(X).
\end{eq}%
So, $\si_0(X)=1$, $\si_1(X)=\tr(X)$ and $\si_n(X)=\det(X)$. The following result holds (see~\cite{Zubkov99}). 

\begin{theo}\label{theo_matrix}
If $\Char{\FF}\neq2$, then the algebra of matrix invariants $\FF[V]^{O(n)}$ is generated by $\si_t(B)$ ($1\leq t\leq n$), where $B$ ranges over all monomials in $X_1,\ldots,X_d$, $X_1^T,\ldots,X_d^T$. Moreover, we can assume that $B$ is primitive, i.e., is not equal to a power of a shorter monomial. 
\end{theo}

\remark\label{rem_tr} In the case of a characteristic zero field it is enough to take
traces instead of $\si_t$, $1\leq t\leq n$, in the formulation of
Theorem~\ref{theo_matrix}.

\bigskip
 
In characteristic zero case Procesi~\cite{Procesi76} described relations between generators for $\FF[V]^G$ for $G\in\{GL(n),O(n),\Sp(n)\}$. Zubkov~\cite{Zubkov96} described relations  for matrix $GL(n)$-invariants over a field of arbitrary characteristic. We have described relations for matrix $O(n)$-invariants (see Theorem~\ref{theo_relations}) in the case of characteristic different from two. The proofs of Theorem~\ref{theo_relations} and Lemma~\ref{lemma1_definition} will be given in a separate paper.

\section{Relations}\label{section_relations}

For a vector $\un{t}=(t_1,\ldots,t_u)\in\NN^u$ we write $\#\un{t}=u$, where ${\NN}$ stands for the set of non-negative integers. In this paper we use the following notions: 
\begin{enumerate}
\item[$\bullet$] the monoid $\M$ (without unity) freely generated by {\it letters}  $x_1,\ldots,x_d,x_1^T,\ldots,x_d^T$, the vector space $\M_{\FF}$ with the basis $\M$, and $\N\subset\M$ the subset of primitive elements, where the notion of a primitive element is defined as above; 

\item[$\bullet$] the involution ${}^T:\M_{\FF}\to\M_{\FF}$ defined by $x^{TT}=x$ for a letter $x$ and $(a_1\cdots a_p)^T=a_p^T\cdots a_1^T$ for $a_1,\ldots,a_p\in\M$;

\item[$\bullet$] the equivalence $y_1\cdots y_p\sim z_1\cdots z_p$ that holds 
if there exists a cyclic permutation $\pi\in S_p$ such that
$y_{\pi(1)}\cdots y_{\pi(p)} = z_1\cdots z_p$ or $y_{\pi(1)}\cdots y_{\pi(p)} = z_p^T\cdots z_1^T$, where $y_1,\ldots,y_p,z_1,\ldots,z_p$ are letters;

\item[$\bullet$] $\M_{\si}$, the ring with unity of (commutative) polynomials over $\FF$ freely generated by the ``symbolic'' elements $\si_t(\al)$, where $t>0$ and $\al\in\M_{\FF}$;   

\item[$\bullet$] $\N_{\si}$, a ring with unity of (commutative) polynomials over $\FF$ freely generated by the ``symbolic'' elements $\si_t(\al)$, where $t>0$ and $\al\in\N$ ranges over $\sim$-equivalence classes; note that $\N_{\si}\simeq \M_{\si}/L$, where the ideal $L$ is described in Lemma~\ref{lemma1_definition} (see below).   




\end{enumerate} 

We will use the following notation:
$$\tr(\al)=\si_1(\al)$$
for any $\al\in\M_{\FF}$.   

For a letter $\beta\in\M$ define
$$X_{\beta}=
\left\{
\begin{array}{rl}
X_{i},&\text{if } \beta=x_i\\
X_{i}^T,&\text{if } \beta=x_i^T\\
\end{array}
\right..
$$
Given $\al=\al_1\cdots \al_p\in\M$, where $\al_i$ is a letter, we assume that $X_{\al}=X_{\al_1}\cdots X_{\al_p}$. Consider a surjective homomorphism of algebras
$$\Psi_n:\N_{\si} \to \FF[V]^{O(n)}$$ %
defined by $\si_t(\al) \to \si_t(X_{\al})$, if $t\leq n$, and $\si_t(\al) \to 0$ otherwise. Its kernel $K_{n}$ is the ideal of {\it relations} for $\FF[V]^{O(n)}$. 
Elements of $\bigcap_{i>0} K_{i}$ are called {\it free} relations.  
For $\al,\be,\ga\in\M_{\FF}$ and $t,r\in\NN$, an element $\si_{t,r}(\al,\be,\ga)\in\N_{\si}$ was introduced in~\cite{ZubkovII} (see Definition~\ref{def1_definition} below).

\begin{theo}\label{theo_relations} If $\Char{\FF}\neq2$, then the ideal of relations $K_{n}$ for $\FF[V]^{O(n)}\simeq  \N_{\si}/K_{n}$ is generated by $\si_{t,r}(\al,\be,\ga)$, where $t+2r>n$ and $\al,\be,\ga\in\M_{\FF}$.
\end{theo}

\bigskip

The following result was known in characteristic zero case.

\begin{cor} In case $G=O(n)$ there are no non-zero free relations.
\end{cor}

\section{The definition of $\si_{t,r}$}\label{section_definition}
In this section we assume that $\algA$ is a commutative unitary algebra over the field $\FF$ and all matrices are considered over $\algA$. 

Let us recall some formulas. In what follows $A,A_1,\ldots,A_p$ stand for $n\times n$ matrices and $1\leq t\leq n$. Amitsur's formula states~\cite{Amitsur_1980}:
\begin{eq}\label{eq_Amitsur}
\si_t(A_1+\cdots+A_p)=\sum (-1)^{t-(j_1+\cdots+j_q)} \si_{j_1}(\ga_1)\cdots\si_{j_q}(\ga_q),
\end{eq}%
where the sum ranges over all pairwise different primitive cycles $\ga_1,\ldots,\ga_q$ in
letters $A_1,\ldots,A_p$ and positive integers $j_1,\ldots,j_q$ with
$\sum_{i=1}^{q}j_i\deg{\ga_i}=t$. Denote the right hand side of~\Ref{eq_Amitsur} by $F_{t,p}(A_1,\ldots,A_p)$. As an example,
\begin{eq}\label{eq_si2}
\si_2(A_1+A_2)=\si_2(A_1)+\si_2(A_2)+\tr(A_1)\tr(A_2)-\tr(A_1A_2).
\end{eq}%
Note that for $a\in\algA$ we have 
\begin{eq}\label{eq20_definition}
\si_t(a A)=a^t\si_t(A).
\end{eq}%

For $l\geq2$ we have the following well-known formula:
\begin{eq}\label{eq_D}
\si_t(A^l)=\sum\limits_{i_1,\ldots,i_{t l}\geq0}b^{(t,l)}_{i_1,\ldots,i_{t l}}
    \si_1(A)^{i_1}\cdots\si_{t l}(A)^{i_{t l}},%
\end{eq}%
\noindent where we assume that $\si_i(A)=0$ for $i>n$. 
Denote the right hand side of~\Ref{eq_D} by $P_{t,l}(A)$. 
In~\Ref{eq_D} coefficients $b^{(t,l)}_{i_1,\ldots,i_{rl}} \in \ZZ$ do not depend on $A$
and $n$. If we take a diagonal matrix $A=\diag(a_1,\ldots,a_n)$, then $\si_t(A^l)$ is
a symmetric polynomial in $a_1,\ldots,a_n$ and $\si_i(A)$ is the
$i^{\rm th}$ elementary symmetric polynomial in $a_1,\ldots,a_n$, where $1\leq i\leq n$. Thus the coefficients
$b^{(t,l)}_{i_1,\ldots,i_{tl}}$ with $tl\leq n$ can easily be found. As an example,
\begin{eq}\label{eq_tr_a2}
\tr(A^2)=\tr(A)^2-2\si_2(A).
\end{eq}%

\begin{lemma}\label{lemma1_definition} We have $\N_{\si}\simeq \M_{\si}/ L$ for the ideal $L$ generated by
\begin{enumerate}
\item[(a)] $\si_t(\al_1+\cdots+\al_p)-F_{t,p}(\al_1,\ldots,\al_p)$, 

\item[(b)] $\si_t(a\,\al)-a^t\si_t(\al)$,

\item[(c)] $\si_t(\al^l)-P_{t,l}(\al)$,

\item[(d)] $\si_t(\al\be)-\si_t(\be\al)$, 
 
\item[(e)] $\si_t(\al)-\si_t(\al^T)$,
\end{enumerate}
where $p>1$, $\al,\al_1,\ldots,\al_p\in \M_{\FF}$, $a\in\FF$, $t>0$, and $l>1$. 
\end{lemma}

\bigskip

A {\it quiver} $\Q=(\Q_0,\Q_1)$ is a finite oriented graph, where $\Q_0$ is the set of
vertices and $\Q_1$ is the set of arrows.  Multiple arrows and loops in $Q$ are allowed.
For an arrow $\al$, denote by $\al'$ its head
and by $\al''$ its tail, i.e., %
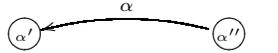
$$\vcenter{
\xymatrix@C=1cm@R=1cm{ %
\vtx{\al'}\ar@/^/@{<-}[rr]^{\al} && \vtx{\al''}\\
}} \quad.
$$
We say that $\al=\al_1\cdots \al_s$ is a {\it path} in $\Q$ (where $\al_1,\ldots,\al_s\in
\Q_1$), if $\al_1''=\al_2',\ldots,\al_{s-1}''=\al_s'$, i.e.,
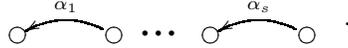
$$\vcenter{
\xymatrix@C=1cm@R=1cm{ %
\vtx{\;}\ar@/^/@{<-}[r]^{\al_1} & %
\vtx{\;} & %
\vtx{\;}\ar@/^/@{<-}[r]^{\al_s} & %
\vtx{\;}\\
}} \quad.
\begin{picture}(0,0)
\put(-73,-3){
\put(0,0){\circle*{2}} %
\put(-5,0){\circle*{2}} %
\put(5,0){\circle*{2}} %
} %
\end{picture}
$$
The head of the path $\al$ is $\al'=\al_1'$ and the tail is $\al''=\al_s''$.  A path
$\al$ is called {\it closed} if $\al'=\al''$.

\definition{of a mixed quiver} A quiver $\Q$ is called {\it mixed} if there are 
two maps ${}^T:\Q_0\to\Q_0$ and ${}^T:\Q_1\to\Q_1$ satisfying  
\begin{enumerate}
\item[$\bullet$] $v^{TT}=v$, $\be^{TT}=\be$;

\item[$\bullet$] $(\be^T)'=(\be'')^T$, $(\be^T)''=(\be')^T$
\end{enumerate} 
for all $v\in\Q_0$ and $\be\in\Q_1$.
\bigskip 

Assume that $\Q$ is a mixed quiver. Denote by $\M(\Q)$ the set of all closed paths in $\Q$ and denote by $\N(\Q)\subset \M(\Q)$ the subset of primitive paths.  Given a path $\al$ in $\Q$, we define the path $\al^{T}$ and introduce $\sim$-equivalence on $\M(\Q)$ in the same way as in Section~\ref{section_relations}. Moreover, we define $\M_{\FF}(\Q)$, $\M_{\si}(\Q)$, and $\N_{\si}(\Q)$ in the same way as $\M_{\FF}$, $\M_{\si}$, and $\N_{\si}$ have been defined in Section~\ref{section_relations}.   

\smallskip
Let $t,r\in\NN$. In order to define $\si_{t,r}$, we consider the quiver $\Q$
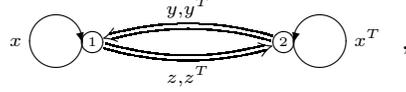
$$
\loopR{0}{0}{x} %
\xymatrix@C=1cm@R=1cm{ %
\vtx{1}\ar@2@/^/@{<-}[rr]^{y,y^T} &&\vtx{2}\ar@2@/^/@{<-}[ll]^{z,z^T}\\
}%
\loopL{0}{0}{x^T}\qquad,
$$
where there are two arrows from vertex $1$ to vertex $2$ and there are two arrows in the
opposite direction. We define $1^T=2$ for vertex $1$ to turn $\Q$ into a mixed quiver.
Denote the degree of $\al\in\M(\Q)$ by $\deg{\al}$, the degree of $\al$ in a letter $\be$ by $\deg_{\be}{\al}$, and the multidegree of $\al$ by 
$$\mdeg(\al)=%
(\deg_{x}{\al}+\deg_{x^T}{\al},\deg_{y}{\al}+\deg_{y^T}{\al},\deg_{z}{\al}+\deg_{z^T}{\al}).$$

\definition{of $\si_{t,r}(x,y,z)$}\label{def1_definition} Denote by  ${\I}={\I}_{t,r}$ the set of pairs $(\un{j},\un{\al})$ such that 
\begin{enumerate}
 \item[$\bullet$] $\#\un{j}=\#\un{\al}=p$ for some $p$;

 \item[$\bullet$] $\al_1,\ldots,\al_p\in \N(\Q)$ belong to pairwise different $\sim$-equivalence classes and $j_1,\ldots,j_p\geq1$;

 \item[$\bullet$] $j_1\mdeg(\al_1)+\cdots+j_p\mdeg(\al_p)=(t,r,r)$.
\end{enumerate}
Then we define $\si_{t,r}(x,y,z)\in\N_{\si}(\Q)$ by
\begin{eq}\label{eq1_definition}
\si_{t,r}(x,y,z)=\sum_{(\un{j},\un{\al})\in{\I}} %
(-1)^{\xi} \;\si_{j_1}(\al_1)\cdots\si_{j_p}(\al_p), %
\end{eq}%
where $p=\#\un{j}=\#\un{\al}$ and $\xi=\xi_{\un{j},\un{\al}}=t+\sum_{i=1}^p j_i(\deg_y{\al_i}+\deg_z{\al_i}+1)$. 
For $t=r=0$ we define $\si_{0,0}(x,y,z)=1$. Moreover, if $\al,\be,\ga\in\M_{\FF}$, then we define $\si_{t,r}(\al,\be,\ga)\in\N_{\si}$ as the result of substitution $x\to\al$, $y\to\be$, $z\to\ga$ in~\Ref{eq1_definition}.

\example\label{ex_DP} {\bf 1.} If $t=0$ and $r=1$, then $\sim$-equivalence classes on $\N(\Q)$ are $y z,\, y z^T,\, \ldots$ %
Hence, $\si_{0,1}(x,y,z)=-\tr(yz)+\tr(yz^T)$.
\smallskip

{\bf 2.} If $t=r=1$, then $\sim$-equivalence classes on $\N(\Q)$ are 
$$x,\, y z,\, y z^T,\, x y z,\, x y z^T,\, x y^T z,\, x y^T z^T,\, \ldots$$
and we can see that $\si_{1,1}(x,y,z)=$
$$-\tr(x)\tr(yz)+\tr(x)\tr(yz^T)+ \\%
\tr(x y z)-\tr(x y z^T)-\tr(x y^T z)+\tr(x y^T z^T).$$ %

\remark\label{rem01} $\si_{t,0}(x,y,z)=\si_t(x)$.
\bigskip

The decomposition formula from~\cite{Lopatin_bplp} implies that for $n\times n$ matrices $A_i$, $i=1,2,3$, with $n=t_0+2r$,  $t_0\geq0$ we have  
\begin{eq}\label{eq1_history}
\DP_{r,r}(A_1+\la E,A_2,A_3)=\sum_{t=0}^{t_0} \la^{t_0-t}\si_{t,r}(A_1,A_2,A_3),
\end{eq}%
where $\DP_{r,r}(A_1,A_2,A_3)$ stands for the determinant-pfaffian (see~\cite{LZ1}) and $\si_{t,r}(A_1,A_2,A_3)$ is defined as the result of the substitution $a_i\to A_i$, $a_i^T\to A_i^T$ in $\si_{t,r}(a_1,a_2,a_3)$. Thus $\DP_{r,r}$ relates to $\si_{t,r}$ in the same way as the determinant relates to $\si_t$.

\section{Application}\label{section_application}

Given an $\NN$-graded algebra $\algA$, denote by $\algA^{+}$ the subalgebra generated by elements of $\algA$ of positive degree. It is easy to see that a set $\{a_i\} \subseteq \algA$ is a minimal (by inclusion) homogeneous system of generators (m.h.s.g.) for $\algA$~if and only if $\{\ov{a_i}\}$ is a basis for $\ov{\algA}={\algA}/{(\algA^{+})^2}$ and $\{a_i\}$ are homogeneous. Let us recall that an element $a\in\algA$ is called {\it
decomposable} if it belongs to the ideal $(\algA^{+})^2$. 
Therefore the
least upper bound for the degrees of elements of a m.h.s.g.~for $\FF[V]^{O(n)}$ is equal to
the maximal degree of indecomposable invariants and we denote it by $D_{\rm max}$. 

As an application of Theorem~\ref{theo_relations} we obtained the following result in~\cite{Lopatin_O3}. 

\begin{theo}\label{theo_D} Let $n=3$ and $d\geq1$. Then
\begin{enumerate} 
\item[$\bullet$] If $\Char\FF =3$, then $2d+4\leq D_{\rm max}\leq 2d+7$.

\item[$\bullet$] If $\Char\FF \neq 2,3$, then $D_{\rm max}=6$.
\end{enumerate}
\end{theo}

Moreover, in~\cite{Lopatin_skewO3} we described a m.h.s.g.~for orthogonal invariants of  $d$-tuples of $3\times 3$ skew-symmetric matrices. As about matrix $GL(n)$-invariants in case $n=3$, its minimal system of generators was explicitly calculated in~\cite{Lopatin_Comm1} and~\cite{Lopatin_Comm2}.


\bigskip
\noindent{\bf Acknowledgements.} The paper has been supported by RFFI 10-01-00383a.

\end{document}